\documentclass[a4paper,12pt]{article}
\usepackage[english]{babel}
\usepackage{amsmath}
\usepackage{amssymb}
\usepackage{tabls}
\usepackage[pdftex]{graphicx}
\usepackage{calc}
\usepackage[T1]{fontenc}
\usepackage[latin1]{inputenc}
\title{Introducing Crystalline Graded Algebras}
\author{Erna Nauwelaerts \\ University of Hasselt \\ \texttt{erna.nauwelaerts@uhasselt.be} \and Freddy Van Oystaeyen \\ University of Antwerp \\ \texttt{fred.vanoystaeyen@ua.ac.be}}
\date{}

\newcommand{\modu}{\textup{ mod }}

\newcommand{\aut}{\textup{Aut}}
\newcommand{\out}{\textup{Out}}
\newcommand{\inn}{\textup{Inn}}
\newcommand{\N}{\mathbb{N}}
\newcommand{\Z}{\mathbb{Z}}

\newcommand{\C}{\mathbb{C}}

\newcommand{\blok}{\hfill \Box}
\newcommand{\str}{\ |\ }

\newtheorem{defi}{Definition}[section]
\newtheorem{opm}[defi]{Remark}
\newtheorem{prop}[defi]{Proposition}
\newtheorem{gev}[defi]{Corollary}

\newtheorem{lem}[defi]{Lemma}

\begin{document}
\maketitle

\begin{abstract}
	We introduce a new class of graded rings extending the class of generalized Weyl algebras. These rings are orders in crossed products of the most general type, and we introduce their basic structure theory. We provide an extensive list of examples, some completely new but also some considered earlier, and highlight some specific structure results. Many new and interesting problems about the crystalline graded rings may be identified.
\end{abstract}

\section*{Introduction}

The theory of group actions on algebras and related crossed product constructions has roots in representation theory of finite groups and the structure theory of central simple algebras or the Brauer group of a field.  Inspired by the success of an approach via integral ring extensions in Number Theory and Local Field Theory a noncommutative theory of orders and maximal orders in central simple algebras was developed.  In this theory orders in crossed products fit!  Modular representation theory e.g. integral representations may be seen as a succesful branch of that same idea.  On the other hand a second important line of application is found in the theory of projective representations, their relation with twisted group rings, Schur multipliers and ray classes providing a link to crystallography explaining the terminology we have chosen.\\
From a different point of view, dealing with generalizations of strongly graded rings (cfr. \cite{NaVO1}, \cite{NaVO2} ), rather unexpected gradations appeared on algebras coming from other mathematical theories, e.g. Weyl algebras, certain quantum groups.  A special case of these generalized crossed products was extensively studied by V. Bavula (cfr. \cite{B1}, \cite{B2} ) who identified these algebras as 'generalized Weyl algebras' (GWA).  Then V. Bavula and the second author exploited fully the graded ring structure of these GWA, in particular the fact that they are graded orders in a twisted Laurent series ring, to arrive at deeper structural results leading to, for example, a classification of simple modules for the second Weyl algebra.\\
In this paper we introduce a general class of graded rings containing as special examples the Weyl algebras and quantizations thereof, the generalized Weyl algebras including 'quantum $\mathfrak{sl}_2$', generalizations of Clifford algebras as well as some specific new interesting new examples.  We view this paper as an introduction of a subject raising a lot of new interesting questions and allowing several new developments in different directions e.g. depending wether the grading group is finite, infinite abelian, totally ordered,...  In the first part we provide the general theory of crystalline graded rings, establishing how they may be viewed as orders in crossed products, the second part contains a lot of examples old and new.\\\\

Both authors supported by a Research Project of the FWO (nr. G.0622.06); second author acknowledges support of the E.C. project Liegrits (RTN 505078).

\section{Definition and basic properties of crystalline graded rings}

Throughout $A$ is an associative ring with unit element $1$, $A_0$ will be a subring of $A$ containing $1$.  Unless otherwise mentioned $G$ will be a completely arbitrary group.  Let $u:G \rightarrow A: g \mapsto u_g$ be a map of sets such that $u_e = 1$ where $e$ is the neutral element of $G$, and $u_g \neq 0$ for all $g \in G$.  We make the following assumptions:

\begin{itemize}
	\item (C1) $A= \mathop \oplus \limits_{g \in G} A_0 u_g$
	\item (C2) For every $g \in G$, $A_0 u_g = u_g A_0$ and this is a free left $A_0$-module of rank one.
	\item (C3) The decomposition in C1 makes $A$ into a $G$-graded ring with $A_0 = A_e$.
\end{itemize}

As a direct consequence of these assumptions we have the following basic lemma.

{\lem \label{1} With conventions and notation as above:
\begin{enumerate}
	\item For every $g \in G$, there is a set map $\sigma_g : A_0 \rightarrow A_0$ defined by: $u_g a_0 = \sigma_g(a_0)u_g$ for $a_0 \in A_0$.  The map $\sigma_g$ is in fact a surjective ring morphism.  Moreover, $\sigma_e = I_{A_0}$.
	\item There is a set map $\alpha : G \times G \rightarrow A_0$ defined by $u_g u_h = \alpha(g,h)u_{gh}$ for $g,h \in G$.  For any triple $g,h,t \in G$ the following equalities hold:
	\begin{eqnarray}
	\alpha(g,h)\alpha(gh,t)&=&\sigma_g(\alpha(h,t))\alpha(g,ht) \label{3}\\
	\sigma_g(\sigma_h(a_0))\alpha(g,h)&=& \alpha(g,h)\sigma_{gh}(a_0) \label{4}
	\end{eqnarray}
\item For all $g \in G$ we have the equalities $\alpha(g,e) = \alpha(e,g) = 1$ and $\alpha(g,g^{-1}) = \sigma_g(\alpha(g^{-1},g)).$
\end{enumerate}
}
\textbf{Proof}
\begin{enumerate}
	\item The definition of $\sigma_g$ as a map of sets follows from (C2) above.  For $a_0, b_0 \in A_0$ we obtain:
	\begin{eqnarray*}
	u_g(a_0 b_0) = (u_g a_0)b_0 &\Rightarrow& \sigma_g(a_0 b_0) = \sigma_g(a_0) \sigma_g(b_0) \\
	u_g(a_0 + b_0) = u_g a_0 + u_g b_0 &\Rightarrow& \sigma_g(a_0 + b_0) = \sigma_g(a_0)+ \sigma_g(b_0)
	\end{eqnarray*}
	using that $A_0u_g$ is a free left $A_0$-module with basis $u_g$.  Observe also that $\sigma_e = I_{A_0}$ because for $a_0 \in A_0: a_0 = u_e a_0 = \sigma_e(a_0)u_e = \sigma(a_0)$.  The surjectivity of $\sigma_g$ follows from $A_0 u_g \subset u_g A_0$, i.e. for an arbitrary $b_0 \in A_0$ we have that $b_0 u_g = u_g a_0$ for some $a_0 \in A_0$, thus $b_0 = \sigma_g(a_0)$.  The surjectivity of $\sigma_g$ implies that $\sigma_g(1) = 1$.
	\item The definition of $\alpha$ as a map of sets follows from (C3).  For any triple $g,h,t \in G$ we find $(u_g u_h)u_t = u_g(u_h u_t)$.  Hence,
	\[\alpha(g,h)\alpha(gh,t)u_{ght} = \sigma_g(\alpha(h,t))\alpha(g,ht)u_{ght}\]
	and the claim follows from the fact that $A_0 u_{ght}$ is a free left $A_0$-module with basis $u_{ght}$.  Secondly, for $a_0$ we have $u_g(u_h a_0) = (u_g u_h)a_0$ and this yields:
	\[\sigma_g(\sigma_h(a_0))\alpha(g,h)u_{gh}=\alpha(g,h)\sigma_{gh}(a_0)u_{gh}\]
	Thus proving the second claim.
	\item First, $u_g = u_g u_e = \alpha(g,e)u_g$ gives $\alpha(g,e) = 1$.  Similarly $\alpha(e,g)=1$.  Furthermore, from
	\[\alpha(g,g^{-1})\alpha(e,g)= \sigma_g(\alpha(g^{-1}, g))\alpha(g,e)\] 
	we find $\alpha(g,g^{-1}) = \sigma_g(\alpha(g^{-1}, g))$. $\blok$
\end{enumerate}

{\gev \label{2}Keep the above conventions and notations and let $g,h \in G$ then:
\[\alpha(g,g^{-1})= \sigma_g(\alpha(g^{-1}, gh))\alpha(g,h)\]}\\
\textbf{Proof}\\
Apply the first equality in Lemma \ref{1}(2) to the triple $g, g^{-1},gh$ and use that $\alpha(e,gh)=1$. $\blok$\\

Let $S(G)$ be the multiplicative set in $A_0$ generated by $\{\alpha(g,g^{-1})\str g \in G\}$ and let $S(G \times G)$ stand for the multiplicative set generated by $\{\alpha(g,h)\str g,h \in G\}$.  The most interesting situation appears when $A_0$ has no $S(G)$-torsion.  However, for the moment we will make a much weaker assumption: $0 \notin S(G)$.  It follows directly from Corollary \ref{2} that $0 \notin S(G)$ implies $0 \notin S(G \times G)$.

{\prop \label{5} With conventions and notations as above:
\begin{enumerate}
	\item For all $g, h \in G$ the $\alpha(g,h)$ are normalizing elements of $A_0$ in the sense that $A_0 \alpha(g,h) = \alpha(g,h) A_0$.
	\item Assume that $0 \notin S(G)$.  Then the multiplicative set $S(G)$ is a left Ore set of $A_0$ and every element of $S(G \times G)$ is invertible in the ring $S(G)^{-1} A_0$.
\end{enumerate}
}
\textbf{Proof}
\begin{enumerate}
	\item Let $g,h \in G$ and look at the equality \ref{4} in Lemma \ref{1}(2), i.e.
	\[\sigma_g(\sigma_h(x_0))\alpha(g,h)= \alpha(g,h)\sigma_{gh}(x_0)\]
	with $x_0 \in A_0$ and $\sigma_g \circ \sigma_h$ as well as $\sigma_{gh}$ being surjective ring morphisms.  For a given $a_0 \in A_0$ we may pick $a'_0 \in A_0$ such that $\sigma_{gh}(a'_0)=a_0$.  Then we obtain $\alpha(g,h)a_0 = \sigma_g(\sigma_h(a'_0))\alpha(g,h)$.  Furthermore, given $b_0 \in A_0$ we pick $b'_0 \in A_0$ such that $\sigma_g(\sigma_h(b'_0)) = b_0$.  Then $b_0 \alpha(g,h) = \alpha(g,h) \sigma_g(b'_0)$.  Observe that for $a_0, b_0 \in A_0$ we may select $a'_0b'_0$ for the representative of $a_0 b_0$ and $\sigma_g\sigma_h(a'_0 b'_0) = \sigma_g \sigma_h(a'_0)\sigma_g \sigma_h(b'_0)$ yields $\alpha(g,h) a_0 b_0 = \sigma_g \sigma_h(a'_0 b'_0) \alpha(g,h)$ holds.
	\item Let us check the left Ore condition for $S(G)$.  First, let $s \in S(G)$ and $a_0 \in A_0$.  Since $s$ is a product of normalizing elements of $A_0$, it is itself a normalizing element of $A_0$.  Therefore there exists an element $a'_0 \in A_0$ such that $sa_0 = a'_0s$.\\
	Secondly, assume that $a_0 s = 0$ for some $s \in S(G)$ and $a_0 \in A_0$.  We have to establish that $s' a_0 = 0$ for some $s' \in S(G)$.  First consider $a_0 \alpha(g,g^{-1}) = 0$ for some $g \in G$; iteration will lead to the case of an arbitrary $s \in S(G)$.  Pick $x_0 \in A_0$ such that $a_0 = \sigma_g(\sigma_{g^{-1}}(x_0))$.  Then we obtain from the second equality in Lemma \ref{1}(2):
	\[a_0 \alpha(g, g^{-1}) = \alpha(g,g^{-1})\sigma_e(x_0)= \alpha(g,g^{-1})x_0\]
Consequently $\alpha(g,g^{-1})x_0 = 0$.  This yields:
\begin{eqnarray*}
0 &=& \sigma_g \circ \sigma_{g^{-1}}(\alpha(g,g^{-1})x_0)=\sigma_g \circ \sigma_{g^{-1}}(\alpha(g,g^{-1}))\sigma_g \circ \sigma_{g^{-1}}(x_0)\\
\ &=& \sigma_g \circ \sigma_{g^{-1}}(\alpha(g,g^{-1})) a_0
\end{eqnarray*}
Applying Lemma \ref{1}(3) we have:
\[\alpha(g,g^{-1}) = \sigma_g(\alpha(g^{-1},g)) = \sigma_g(\sigma_{g^{-1}}(\alpha(g,g^{-1})))\]
So we obtain $\alpha(g,g^{-1})a_0 = 0)$.\\
Now we write $s = \alpha(g_1, g_1^{-1})\ldots\alpha(g_n, g_n^{-1})$, where repetition of $g_i$ is allowed, and we assume that $a_0 s = 0$.  Application of the foregoing leads to:
\[\alpha(g_n, g_n^{-1})a_0\alpha(g_1, g_1^{-1})\ldots\alpha(g_{n-1}, g_{n-1}^{-1})=0\]
Repetition gives $\alpha(g_1, g_1^{-1})\ldots\alpha(g_n, g_n^{-1})a_0 = 0$, and thus $sa_0 = 0$.\\
Let $g \in G$.  Corollary \ref{2} yields: if $\alpha(g,g^{-1})$ is invertible in some overring (of an image) of $A_0$, then so is $\alpha(g,h)$ for all $h \in G$.  Indeed, from $\alpha(g,g^{-1}) = \sigma_g(\alpha(g^{-1}, gh))\alpha(g,h)$ it follows that $\alpha(g,h)$ is left invertible in $S(G)^{-1}A_0$ say $t \alpha(g,h) = 1$.  There is an $s \in S(G)$ such that $st \in A_0$ and $a' \in A_0$ such that $st \alpha(g,h) = \alpha(g,h)a'$ because $\alpha(g,h)$ is normalizing.  So $\alpha(g,h)a'=s$ but then $\alpha(g,h)a's^{-1} = 1$ with $a's^{-1} \in S(G)^{-1}A_0$. $\blok$
\end{enumerate}

{\gev{\textbf{(of the proof)}} For every $g \in G$ we have that $\ker(\sigma_g \sigma_{g^{-1}}) = t_{\alpha(g, g^{-1})}(A_0)$, the $\alpha(g,g^{-1})$-torsion part of $A_0$.}\\
\textbf{Proof}\\
Since $\alpha(g,g^{-1})$ is normalizing, that set of elements of $A_0$ left annihilated by $\alpha(g,g^{-1})$ is an ideal of $A_0$.  More generally, $\sigma_g \sigma_h (z)=0$ for $z \in A_0$ is equivalent to $\alpha(g,h)\sigma_{gh}=0$, hence $t_{\alpha(g, g^{-1})}(A_0)$ is exactly $\sigma_{gh}(\ker \sigma_g \sigma_h)$.  Observe that the latter would be zero in case $\sigma$ is a group morphism to $\textup{End}(A_0)$.  $\blok$\\

The ring extension $A$ of $A_0$ is now not too bad since it is a free left $A_0$-module generated by normalizing elements (some properties of such extensions are known from work of L. Small and C. Robson on so-called liberal extensions, cfr. \cite{RS}).

{\prop With notation and conventions as before: $S(G)$ is a left Ore set of $A$.}\\
\textbf{Proof}\\
For $g, h \in G$ we have established earlier from Lemma \ref{1}(2) (\ref{3}) with $t = h^{-1}$:
	\[\alpha(g,h)\alpha(gh,h^{-1})=\sigma_g(\alpha(h,h^{-1}))\]
We see that for all $g \in G$, the $\sigma_g$-orbit of elements of $S(G)$ stays within $S(G \times G)$, while moreover from Proposition \ref{5}(2) elements of $S(G \times G)$ are invertible in $S(G)^{-1}A_0$, consequently for every $d \in S(G \times G)$ there exists an $s \in S(G) \cap A_0 d$.  This holds in particular for $d$ that is a product of $\sigma_g$-images of elements of $S(G)$.  Now look at $a \in A$ and $s \in S(G)$ such that $as=0$.  There is a unique homogeneous decomposition
\[a = \sum_{g \in G}^{} {}^{'} a_g u_g\]
with $a_g \in A_0$.  From $as = 0$ then follows that $a_g u_g s = 0$ for all $g \in G$.  We have $a_g u_g = u_g a'_g$ for some $a'_g \in A_0$ and $u_g a'_g s = 0$ yields $\sigma_g(a'_g s)=0$, hence $\sigma(a'_g)\sigma_g(s)=0$.  Now in $S(G)^{-1}A_0$, all elements $\sigma_g(s)$ are invertible as we observed above, hence the $\sigma(a'_g)$ are $S(G)$-torsion elements of $A_0$, say $s_g \sigma(a'_g)=0$ for suitable $s_g \in S(G)$, $g \in G$.  Since only finitely many $a_g$ appeared in the homogeneous decomposition of $a$, only finitely many $s_g$ appear here and thus we may select an $s'' \in S(G)$ such that $s'' \sigma_g(s'_g) = 0$ for all $g$ appearing in the decomposition of $a$.  Hence $s''\sigma_g(a'_g) u_g = 0$ or $s''u_g a'_g = s'' a_g u_g = 0$ for all $g$ as before.  Consequently $s''a = 0$ with $s'' \in S(G)$ as desired.\\
For the other Ore condition, consider $x \in A$, $s \in S(G)$, say
\[x = \sum_{g \in G}^{} {}^{'} u_g x_g = \sum_{g \in G}^{} {}^{'}\sigma_g(x_g)u_g\]
For each $g$ such that $x_g \neq 0$ pick $s_g \in S(G)$ such that $s_g x_g = x'_g s$ (possible, since $S(G)$ is a left Ore set in $A_0$) with $x'_g \in A_0$.  Observe:
\[\sigma_g(s_g)u_g x_g = u_g s_g x_g= u_g x'_g s \in As\]
We have pointed out earlier that there is an element in $S(G) \cap A_0\sigma_g(s_g)$ and since we are looking at only finitely many nonzero elements there also is an $s' \in S(G)$ in $\bigcap_g A_0 \sigma_g(s_g)$.  We then calculate that $s'u_gx_g \in As$, for all $g$ such that $x_g \neq 0$, consequently $s'x \in As$ and the second left Ore condition follows.$\ \ \ $ $\blok$

{\opm{\ }}
\ 
\begin{enumerate}
	\item Our assumption $0 \notin S(G)$ means in particular that $\alpha(g, g^{-1}) \neq 0$ for all $g \in G$.  It follows directly from Corollary \ref{2} that $\alpha(g,h) \neq 0$ for all $g,h \in G$.
	\item If we generate the $\{\sigma_g \str g \in G\}$-invariant multiplicative set in $A_0$ by the images of elements in $S(G)$ then this will be contained in $S(G \times G)$ and therefore consists of elements invertible in $S(G)^{-1}A_0$.  Even if $S(G)$ is not invariant under the $\sigma_g, g \in G$, it behaves 'as if invariant' in view of the forementioned property.  Observe moreover that $S(G \times G)$ need not be a left Ore set in $A$.
\end{enumerate}
	
{\gev The following are equivalent
\begin{enumerate}
	\item $A_0$ is $G(S)$-torsionfree
	\item $A$ is $S(G)$-torsionfree
	\item $\alpha(g,g^{-1})a_0$ for some $g \in G$ implies $a_0 = 0$
	\item $\alpha(g,h)a_0$ for some $g,h \in G$ implies $a_0 = 0$
	\item $A_0 u_g = u_g A_0$ is also free as a right $A_0$-module with basis $u_g$ for every $g \in G$
	\item for every $g \in G$, $\sigma_g$ is bijective hence a ring automorphism of $A_0$
\end{enumerate}
}
	
For a ring monomorphism $A\hookrightarrow B$ we define a subgroup $\aut_B A$ in $\aut A$ consisting of the automorphisms of $A$ that extend to an automorphism of $B$.  An automorphism of $A$ is \textbf{B-inner} if it is induced by an inner automorphism of $B$.  The group of $B$-inner automorphisms of $A$ is a normal subgroup of $\aut_B A$ and we put $\out_B(A)= \aut_B A/\inn_B A$, where $\inn_B A$ stands for the group of $B$-inner automorphisms.

{\gev If $A_0$ is $S(G)$-torsionfree then $\sigma : G \mapsto \aut A_0: g \mapsto \sigma_g$ induces a group homomorphism $\bar \sigma : \out_{S(G)^{-1}A_0}(A_0)$.}\\
\textbf{Proof}\\
Observe that $u_g u_{g^{-1}} = \alpha(g,g^{-1})$ implies that each $u_g$ represents a unit in $S(G)^{-1}A_0$ and for $a_0 \in A_0$ and for $a_0 \in A_0$ we have that $\sigma_g(a_0) = u_g a_0 u_g^{-1}$, hence $\sigma_g \in \aut_{S(G)^{-1}A_0}(A_0)$ for every $g \in G$.  Moreover, $\sigma_g \sigma_h$ and $\sigma_{gh}$ are different just by the inner automorphisms of $S(G)^{-1}A_0$ defined by $\alpha(g,h)$ which is also a unit in that ring. $\blok$

{\gev If $A_0$ is $S(G)$-torsionfree then the following statements are equivalent:
\begin{enumerate}
	\item For $g \in G$, $\alpha(g^{-1}, g) \in Z(A_0)$, resp. $\alpha(g,h) \in Z(A_0)$
	\item For $g \in G$, $\sigma_g \sigma_{g^{-1}}= \sigma_{g^{-1}}\sigma_g = I$, resp $\sigma_g \sigma_h = \sigma_{gh}$
\end{enumerate}}
{\raggedright \textbf{Proof}}\\
For $g \in G$ we have the following sequence of implications:
\begin{eqnarray*}
\sigma_{g^{-1}}\sigma_g = I &\Rightarrow& \alpha(g^{-1},g) \in Z(A_0) \Rightarrow \sigma_g(\alpha(g,g^{-1})) \in Z(A_0)\\
&\Rightarrow& \alpha(g,g^{-1}) \in Z(A_0) \Rightarrow \sigma_g \sigma_{g^{-1}} 
\end{eqnarray*} 
The other claims are obvious from earlier observations, e.g. $\sigma_g \sigma_h(a_0) \alpha(g,h)=\alpha(g,h)\sigma_{gh}(a_0)$ for all $a_0 \in A_0$ as in the equality \ref{4} in Lemma \ref{1}(2).  So $\sigma_g \sigma_h = \sigma_{gh}$ leads to $\alpha(g,h) \in Z(A_0)$ in view of the surjectivity of the maps, conversely $\alpha(g,h) \in Z(A_0)$ leads to $\sigma_g \sigma_h = \sigma_{gh}$. $\blok$\\

If $S(G)$-torsion exists in $A_0$ then we aim to reduce by it so as to arrive in the $S(G)$-torsionfree situation.  Put
\[J = t_{S(G)}(A_0)=\{a_0 \in A_0, sa_0 = 0 \textrm{ for some } s \in S(G)\}\]
Because $S(G)$ is also a left Ore set of $A$ it follows that
\[t_{S(G)}(A)=\{a \in A, sa = 0 \textrm{ for some } s \in S(G)\}\]
is a two sided ideal of $A$.  Clearly, $t_{S(G)}(A)$ is a graded ideal because if $sx=0$ with $x = \sum x_g u_g$, then $sx_g =0$, for all $x_g$, follows and consequently $x_g \in J$ follows too.  From the foregoing we obtain
\[t_{S(G)}(A)=\mathop  \oplus \limits_{g \in G} J u_g = JA\]
We claim that $AJ = JA$.  Indeed, first check $AJ \subset JA$ by looking at $u_g t_0$ with $t_0 \in J$ for some $g \in G$; say $st_0=0$ for some $s \in S(G)$ then $0=u_g s t_0 = \sigma_g(s)u_g t_0$ and since $\sigma_g(s)$ is invertible in $S(G)^{-1}A_0$ it follows that $u_gt_0$ is $S(G)$-torsion in $A$ and thus in $JA$.  Secondly, looking at a homogeneous $t_0 u_g$ with $t_0 \in J$, then there is a $\omega_0 \in A_0$ such that $t_0 u_g = u_g \omega_0$ and $\sigma_g(\omega_0)=t_0$ since $\sigma_g$ is surjective.  Then from $st_0=0$ for some $s \in S(G)$ we derive $st_0u_g = 0 = su_g\omega_0$, hence $0 = u_{g^{-1}}su_g\omega_0 = \sigma_{g^{-1}}(s)\alpha(g^{-1},g)\omega_0$ where $\sigma_{g^{-1}}(s)\alpha(g^{-1},g)$ is invertible in $S(G)^{-1}A_0$, yielding $\omega_0 \in J$ and thus $t_0u_g \in AJ$, finally $JA = AJ$.

{\gev For $t_{S(G)}(A)$ we have that it is generated by $t_{s(G)}(A)_0 = J$ (on the right and on the left) i.e. we have $t_{S(G)}(A) =AJ=JA$ and $J$ is globally $\sigma_g$-invariant for every $g \in G$.  Moreover $\bar A = A/JA$ is graded by $G$
\[\bar A = \mathop  \oplus \limits_{g \in G}(A_0/J)\bar u_g\]
where $\bar u_g = u_g \modu JA$, $\bar A_0 = A_0/J$ amd $(\bar A)_e = \bar A_0$.}\\\\
Observe that an arbitrary graded ideal of $A$, say $I$, need not be generated by $I_0 = I \cap A_0$ either on the left or on the right but we do have that $I/I_0A$ is $G(S)$-torsion.\\
Any $G$-graded ring $A$ with properties C1, C2, C3, and which is $G(S)$-torsionfree is called a \textbf{crystalline graded ring}.  In case $\alpha(g,h) \in Z(A_0)$, or equivalently $\sigma_{gh}=\sigma_g \sigma_h$, for all $g,h \in G$, then we say that $A$ is \textbf{centrally crystalline}.

{\lem If $A$ is $G$-graded with properties C1, C2, C3, then $\bar A = A/JA$ is crystalline graded.}\\
\textbf{Proof}\\
It is clear that $\bar A = \mathop  \oplus \limits_{g \in G} \bar A_0 \bar u_g$ with $\bar u_g \neq 0$ and $\bar u_e = 1$.  For $\bar a_0 \in \bar A_0$ we have for all $g \in G$: $\bar u_g \bar a_0 = \bar \sigma_g(\bar a_0) \bar u_g$ wher $\bar \sigma_g$ is induced by $\sigma_g$ on the quotient $\bar A_0$ (using the $\sigma_g$-invariance of $J$).  That $\bar A_0 \bar u_g$ is free (left) with basis $\bar u_g$ follows because $a_0 u_g \in JA$ if and only if $a_0 \in J$.  Hence $\bar A$ satisfies the same properties as $A$ and is moreover $G(S)$-torsionfree. $\blok$\\

Often the condition $0 \notin S(G)$ or even torsionfreeness with respect to $S(G)$ follow from structural properties available in fairly common situations as in the following situation.

{\prop If $A$ satisfies properties C1, C2 and C3 such that $A_0$ is a prime left Goldie ring then $A$ is crystalline graded.}\\
\textbf{Proof}\\
Since each $\alpha(g,h)$ for $g,h \in G$ is normalizing, $A_0 \alpha(g,h)$ is an ideal of $A_0$ and therefore it contains a regular element of $A_0$, say $z = a_0 \alpha(g,h) = \alpha(g,h)b_0$.  Then $\alpha(g,h)$ is regular too, for every $g,h \in G$, hence $0 \notin S(G)$ and $A_0$ is $S(G)$-torsionfree. $\blok$

{\prop If $A$ satisfies C1, C2 and C3 such that $0 \notin S(G)$ then $S(G)^{-1}A$ is a strongly graded ring.  In this generality $\sigma$ does define a group morphism $G \rightarrow \textup{Pic}(S(G)^{-1}A_0)$ defining an action of $G$ on $Z(S(G)^{-1}A_0)$.}\\
\textbf{Proof}\\
Clearly $\bar A$ is $S(G)$-torsionfree and $S(G)^{-1}A = S(G)^{-1}\bar A$.  Note $S(G)^{-1}\bar A = \mathop  \oplus \limits_{g \in G} (S(G)^{-1}A_0)\bar u_g$ and from $\bar u_g \bar u_{g^{-1}} = \overline{\alpha(g,g^{-1})}$ it follows that $\bar u_g$ is a unit of $S(G)^{-1}A$ contained in $(S(G)^{-1}A)_g$ for all $g \in G$.  Thus
\[(S(G)^{-1}A)_g (S(G)^{-1}A)_h = (S(G)^{-1}A)_{gh}\]
for all $g,h \in G$. $\blok$\\

The structure of $S(G)^{-1}A$ in the foregoing proposition is close to a crossed product, but here we do not have that $\sigma$ defines an action of $G$ on $\bar A_0$ via automorphisms, but it defines a 'projective action' in the sense that:
\[\tau_{gh}\sigma_{gh}=\sigma_g \sigma_h\]
where $\tau_{gh}$ is the morphism associated to $\alpha(g,h)$ induced by an inner of $S(G)^{-1}A_0$.  Observe also that $\sigma$ defines group morphisms:
\[G \rightarrow \textup{Pic}(S(G)^{-1}A_0) \rightarrow \aut(Z(S(G)^{-1}A_0))\] 

\section{Examples of Crystalline Gradations}

\subsection{Generalities}

In most applications we are dealing with $A_0$ that are $k$-algebras for some central (in $A$) field $k$.  Then
\[k[S(G)]=k[\alpha(g,g^{-1}), g \in G]\]
\[k[S(G \times G)]=k[\alpha(g,h), g,h \in G]\]
are subrings of $A_0$ (noncommutative) deserving special attention.  Several classes of crystalline graded $k$-algebras may be defined in connection to forementioned algebras:\\
\textbf{Class 1} \quad For $g,h,t \in G$, $\sigma_t(\alpha(g,h)) \in k[S(G\times G)]$, observe that $\sigma_t(\alpha(g^{-1},g)) \in k[S(G\times G)]$ holds for crystalline graded rings (Lemma \ref{1}(2) (\ref{3}) with $t = h^{-1}$).\\
\textbf{Class 2} \quad For $g,h \in G$, $\alpha(g,h) \in Z(A_0)$, i.e. $A$ is centrally crystalline graded and thus $\sigma_g \sigma_h = \sigma_{gh}$ and moreover we have $\alpha(g,h) \in k[S(G)^{\sigma}]$ where $S(G)^{\sigma} = \{\sigma_h(\alpha(g^{-1},g)), h \in G\}$.\\
\textbf{Class 3} \quad For $g,h \in G$, $\alpha(g,h) \in k[S(G)]$, i.e. $k[S(G \times G)]=k[S(G)]$.  Eventually we may restrict this further to the centrally crystalline graded case.\\

In case $A$ is $S(G)$-torsionfree then certain multiplicative conditions are more natural perhaps.  In $S(G)^{-1}A_0$ we may consider the groups generated by $S(G)$, say $\left\langle S(G)\right\rangle$, resp. by $S(G \times G)$, say $\left\langle S(G \times G)\right\rangle$.  From
\[\alpha(g,h)\alpha(gh,h^{-1})=\sigma_g(\alpha(h,h^{-1}))\]
for $h,g \in G$ it follows that $\left\langle S(G)^\sigma\right\rangle \subset \left\langle S(G \times G)\right\rangle$, from
\[\sigma_g(\alpha(h,t))=\alpha(g,h)\alpha(gh,t)\alpha(g,ht)^{-1}\]
it follows that $\left\langle S(G\times G)^{\sigma}\right\rangle=\left\langle S(G \times G)\right\rangle$ or $\left\langle S(G \times G)\right\rangle$ is $\sigma_g$-invariant for all $g \in G$ (where $\sigma_g$ may be viewed as the canonical extension of $\sigma_g$ on $A_0$ to $S(G)^{-1}A_0$).  Put $C(G) = \left\langle S(G)^{\sigma}\right\rangle \cap A_0$, resp. $C(G \times G) = \left\langle S(G \times G)^{\sigma}\right\rangle \cap A_0$.  Then both $C(G)$ and $C(G \times G)$ are invariant (with respect to $\{\sigma_g, g \in G\}$) multiplicative sets containing $1$ but not $0$.

{\opm Let $A$ be a crystalline graded ring such that $A$ is $S(G)$-torsionfree, then $C(G)$ and $C(G \times G)$ are Ore sets in $A_0$ such that $C(G)^{-1}A_0 = C(G \times G)^{-1}A_0 = S(G)^{-1}A_0$, and both $C(G)$ and $C(G \times G)$ are $\{\sigma_g, g \in G\}$-invariant.}\\
\textbf{Proof}\\
Note that in the $S(G)$-torsionfree case now $S(G \times G)$ is an Ore set; one easily verifies that $\left\langle S(G)^\sigma\right\rangle$ and $\left\langle S(G \times G)\right\rangle$ are Ore sets in $S(G)^{-1}A_0$ and also that $C(G)$ resp. $C(G \times G)$ are invariant Ore sets in $A_0$. $\blok$

\subsection{Examples}
\subsubsection{Crossed products} \label{6}

Given a group $G$, a ring $A_0$ and a group morphism $\sigma:G \rightarrow \aut A_0$.  On $A = \mathop  \oplus \limits_{g \in G} A_0 u_g$, where $A_0 u_g$ is just notation for a copy of $A_0$.  We define for $g \in G, a \in A_0$:
\begin{eqnarray*}
u_g a &=& \sigma_g(a)u_g\\
u_g u_h &=& \alpha(g,h)u_{gh}
\end{eqnarray*}
where $\alpha: G \times G \rightarrow U(A_0)$ the units of $A_0$.  It must satisfy for $g,h,t \in G$:
\[\alpha(g,h)\alpha(gh,t)=\sigma_g(\alpha(h,t))\alpha(g,ht)\]
In this way $A$ is strongly graded by $G$, $A = A_0 * G$ with the ring structure defined as above is called the crossed product $(A_0, G, \sigma, \alpha)$.

\subsubsection{Generalized twisted group rings (see \cite{NVO3})} \label{7}

In \ref{6} we restrict to $\sigma_g$, for all $g \in G$, being the identity, but we allow $\alpha: G \times G \rightarrow A_0 \backslash \{0\}$  Then we proceed as in \ref{6} and denote the obtained (algebra) ring $A$ by $A = A_0^\alpha G$.\\
In classical examples of \ref{6} and \ref{7} there are further restrictions: $A_0 \subset Z(A)$, $\alpha:G \times G \rightarrow U(Z(A_0))$.

\subsubsection{The Weyl algebra $A_1(\C)$} \label{8}

This well-studied algebra may be defined as the quotient of the free algebra $\C\left\langle X,Y\right\rangle$ by the two-sided ideal generated by $YX-XY-1$, it turns out to be isomorphic to a ring of differential operators on a polynomial ring $\C[x]$ generated by two operators i.e. $x$ viewed as multiplication by $x$ in $\C[x]$ and $y$ as $\frac{\partial}{\partial x}$ on $\C[x]$.  Then
\[A_1(\C)=\frac{\C\left\langle X,Y\right\rangle}{(YX-XY-1)} \cong \C[x]\left[y, \frac{\partial}{\partial x}\right]\]
We define a $\Z$-gradation on $A_1(\C)$ by putting $\deg x =1, \deg y =-1$, and therefore:
\begin{eqnarray*}
A_1(\C)_0 &=& \C[xy]\\
A_1(\C)_n &=& \C[xy]x^n, \qquad \textup{for } n \geq 0\\
A_1(\C)_m &=& \C[xy]y^{-m}, \qquad \textup{for } m \leq 0
\end{eqnarray*}
We set $u_n = x^n$ if $n \geq 0$ and $u_m = y^{-m}$ if $m \leq 0$.  We will note $\sigma_{n}$ as $\sigma_{x^n}$ for $n \geq 0$, and $\sigma_{m}$ as $\sigma_{y^{-m}}$ if $m \leq 0$.  It is clear that $\sigma_x(xy)=xy-1$, because $x(xy)=(xy-1)x$, $\sigma_y(xy)=xy+1$, because $y(xy)= (1+xy)y$.  Let us put $t=xy$ for convenience, then we have
\[\sigma:\Z \rightarrow \aut_{\C} \C[t]:n \mapsto (t\mapsto t-n)\]
It is easy to calculate
\begin{multline*}
\alpha(n,-n)=x^ny^n=x^{n-1}ty^{n-1}=(t-n+1)x^{n-1}y^{n-1}=\ldots\\=(t-n+1)(t-n+2)\cdot\ldots\cdot(t)
\end{multline*}
Furthermore  $\alpha(n,-m)$ with $n > m$ ($n,m \in \N$) can be calculated from
\[x^ny^m=x^{n-m}x^my^m=x^{n-m}\alpha(m,-m)=\sigma_x^{n-m}(\alpha(m,-m))x^{n-m}\]
and so
\[\alpha(n,-m)=\sigma_x^{n-m}(\alpha(m,-m))\]
and so on.  Observe that in this example $S(G \times G)$ is the $\sigma$-invariant multiplicative system generated by $S(G)$.

\subsubsection{The quantum Weyl algebra $A_1(q)$} \label{9}

Following argumentation as in \ref{8} but now with respect to the defining relation $YX-qXY-1$.  Again $A_0 = \C[t]$ with $t=xy$ but now $\sigma_{x}:t \mapsto q^{-1}(t-1)$ and so on.  Similar calculations may be carried out leading to formulae in $q$-numbers.

\subsubsection{The quantum plane} \label{10}

The quantum plane is given by generators and relations as
\[A = \frac{K\left\langle X,Y\right\rangle}{XY-\lambda YX}\]
Putting $\alpha(1,-1) = xy =t$, we find that $\sigma_x(t)=\lambda t$ and all calculations are done similarly as above.\\

In many cases we deal with algebras over some base field $k$ and tensor products then provide new examples, e.g. higher Weyl algebras or 'products' of quantum planes.

\subsubsection{Generalized Weyl algebras (cfr. \cite{B1},\cite{B2},\cite{BVO1},\cite{BVO2},\ldots)} \label{11}

This class contains variations on quantum groups like quantum $\mathfrak{sl}_2$ or the Witten algebra, or more generally Rees rings of generalized gauge algebras (cfr. \cite{VO2}).  Even though generalized Weyl algebras are special examples of $\delta$-strongly graded rings with respect to an invertible ideal $\delta$ in degree zero, they were first studied in detail by V. Bavula in a series of papers starting in 1991, e.g. \cite{B1},\cite{B2},\ldots.  Adapting terminology and notation from \cite{BVO1} let us show here how the generalized Weyl algebras become crystalline graded rings.  Consider a ring $D$ (commutative in many examples), and $\sigma = (\sigma_1, \ldots, \sigma_n)$ a set of commuting automorphisms of $D$.  Let $a = (a_1, \ldots, a_n)$ be an $n$-tuple with nonzero entries in $Z(D)$ such that $\sigma_i(a_j)=a_j$ for $i \neq j$.  We define the generalized Weyl algebra (GWA), $A = D(\sigma, a)$, as the ring generated by $D$ and $2n$ symbols $X_1^+, \ldots, X_n^+, X_1^-, \ldots, X_n^-$ satisfying the following rules:
\begin{enumerate}
	\item For $i = 1, \ldots, n : X_i^-X_i^+ = a_i$, $X_i^+ X_i^- = \sigma_i(a_i)$
	\item For all $d \in D$ and all $i$, $X_i^{\pm}d = \sigma_i^{\pm}(d)X_i^{\pm}$
	\item For $i \neq j$, $[X_i^-, X_j^-]=[X_i^+, X_j^+]=[X_i^+, X_j^-]=0$
\end{enumerate}
Now write for $m\in \N$, $v_{\pm}(i) = (X_i^{\pm})^m$.  For $k = (k_1, \ldots, k_n) \in \Z^n$ we set $v_k = v_{k_1}(1)\cdot\ldots\cdot v_{k_n}(n)$.  Putting $A = \mathop \oplus \limits_{k \in \Z^n} A_k$ where $A_k = Dv_k$, we obtain a $\Z^n$-graded ring which we call a GWA of degree $n$.  If $D$ is a $k$-algebra, $k$ invariant under all automorphisms involved, then $A$ is a $k$-GWA.  For GWA's $A$ and $A'$ over $k$, $A \otimes_k A'$ is again a GWA.\\
Let us look at $A = D(\sigma, a)$ of degree $1$, with $a$ nonzero in $Z(D)$, $\sigma \in \aut D$ (this will be a centrally crystalline case).  We may view this special case as an algebra generated over $D$ by $X$ and $Y$ satisfying:
\begin{itemize}
	\item $Xd = \sigma(d)X$, $Yd = \sigma^{-1}(d)Y$ for all $d \in D$
	\item $YX = a$ and $XY = \sigma(a)$
	\item $A = \mathop \oplus \limits_{n \in \Z} A_n$, $A_n = Dv_n$, $v_n = X^n$ for $n \geq 0$ and $v_m = Y^{-m}$ if $m < 0$.  Therefore we obtain $v_n v_m = (n,m) v_{n+m}$.  The following rules define $(n,m)$:
	
\begin{itemize}
	\item For $n > 0$ and $m > 0$ we have $(n,m) = 1$.
	\item Let $n \geq m \geq 0$
	\begin{eqnarray*}
	(n,-m) &=& \sigma^n(a)\cdot\ldots\cdot\sigma^{n-m+1}(a)\\
	(-n,m) &=& \sigma^{-n +1}(a)\cdot\ldots\cdot\sigma^{-n+m}(a)
	\end{eqnarray*}
	\item Let $0 \leq n \leq m$
	\begin{eqnarray*}
	(n,-m) &=& \sigma^n(a)\cdot\ldots\cdot\sigma(a)\\
	(-n,m) &=& \sigma^{-n +1}(a)\cdot\ldots\cdot a
	\end{eqnarray*}
\end{itemize}
\end{itemize}

Let us mention some more specific cases, drawing from \cite{BVO1}.  Take $D = K[t]$, up to a change of variable every $K$-automorphism of $K[t]$ is either $t \mapsto t-1$ or $t \mapsto \lambda t$ for some nonzero $\lambda$ in $K$.  Consider the GWA $A = K[t](\sigma, a)$, $\sigma(t)=t-1$.  Then $A$ is simple if and only if $\textup{char}K = 0$ and there does not exist an irreducible polynomial $p$ in $K[t]$ such that $p$ and $\sigma^{i}(p)$ are multiples of $a$ for some nonzero $i$.\\
Another specific case is obtained by taking $B = K[t,t^{-1}](\sigma,a)$ where $\sigma(t)=\lambda t$ and $\lambda \neq 0, 1$ is in $K$.  Up to changing $X$ to $Xt^i$, $Y$ to $Y$, $t$ to $t$, we may assume that $a \in K[t]$ has a nonzero constant term $a(0)$.  Then $B$ is simple if and only if there does not exist an irreducible $p \in K[t]$ such that $p$ and $\sigma^i(p)$ are multiples of $a$ for some $i$ and $\lambda$ is not a root of unity.  In case $a = 1$ we have that $B$ is localization of the quantum plane $K\left\langle X,Y\right\rangle/(XY-\lambda YX)$ at the multiplicative set generated by $t = XY$.

\subsubsection{Cyclic invariants of the first Weyl algebra} \label{12}

Let $\omega$ be a primitive $m$-th root of unity, let $G$ be the cyclic group of order $m$.  We may let $G$ act on the first Weyl algebra $A_1(K)=K\left\langle x,y\right\rangle$ by $y \mapsto \omega y$, $x \mapsto \omega^{-1}x$, and consider the fixed algebra $A_1(K)^G$.  One easily verifies that
\[A_1(K)^G \cong K[t]\left(\sigma, a = m^m t \left(t + \frac{1}{m}\right)\cdot\ldots\cdot\left(t + \frac{m-1}{m}\right)\right)\]

\subsubsection{Associated graded algebras of some GWA}

Let $A$ be the algebra at the end of \ref{11}, where  $a = \alpha t^k + \textup{(terms of lower degrees)}$.  Then we define a filtration $FA$ on $A$, $F_m A = \sum Kt^iv_j$ with $\kappa(i) +\tau(j) \leq m$, where $\kappa$ and $\tau$ are some weight applied to the indexes.  Then $G_F(A) = K[t][\textup{Id}, \alpha t^k]$, again a GWA.  A simmilar construction, but to be modified carefully, works in case of an automorphism $t\mapsto \lambda t$.

\subsubsection{Quantum $\mathfrak{sl}_2$}

Consider the enveloping algebra $U = U_K(\mathfrak{sl}_2)$ and $c$ the Casimir element of $U$.  Put for $\lambda \in K$
\[U(\lambda)=\frac{U}{U(c-\lambda)}\cong K[t][\sigma, a = \lambda - t(t-1)]\]
with $\sigma(t) = t-1$.  Then $U(\lambda)$ is simple if and only if $\lambda \notin \left\{\frac{n^2 -1}{4}, n =1,2,\ldots \right\}$.\\
The quantized version of $U$ is a well-known quantum group $U_q(\mathfrak{sl}_2)$ where $q \in K$ is not a root of unity.  Then we get
\[U_q(\lambda) = \frac{U_q(\mathfrak{sl}_2)}{(c-\lambda)}\]
with $\lambda \in K$.  These are GWA as follows:
\[U_q(\lambda) \cong K[t,t^{-1}](\sigma, a = \lambda +c)\]
where $\sigma(t)=qt$ and for some $h \in K$
\[c = \left(\frac{t^2}{q^2-1}-\frac{t^{-2}}{q^{-2}-1}\right)(2h)^{-1}\]
We know that $U_q(\lambda)$ is simple if and only if it has no simple finite dimensional module, if and only if for each root of $a$ (take $K = \C$ here), say $\mu$, no $q^i \mu, 0 \neq i \in \Z$ is again a root of $a$.\\
Further examples may be derived from tensor products
\[\frac{U(\mathfrak{sl}_2 \times \ldots \times \mathfrak{sl}_2)}{(c_1-\lambda_1, \ldots, c_n - \lambda_n)}\cong \mathop  \otimes \limits_{i = 1}^n U(\lambda_i)\]

\subsubsection{An example of Bavula, Bekkert as a crystalline ring}

Take $A_0 = K[t]$, $\sigma(t)=t-1$, $a = \alpha(-1,1) = 3^3 t (t-1/3)(t-2/3)$, $\sigma(a)=\alpha(1,-1)=3^3(t-1)(t-4/3)(t-5/3)$.  It is then easy to calculate for $r > 0$:
\begin{eqnarray*}
\sigma^r(a)&=&3^3(t-r)(t-r-1/3)(t-r-2/3)\\
\sigma^{-r}(a) &=& 3^3 (t+r)(t+r-1/3)(t+r-2/3)
\end{eqnarray*}
Furthermore:
\begin{eqnarray*}
\alpha(n,-n) &=& \sigma^n(a)\cdot\ldots\cdot\sigma(a)\\
\alpha(-n,n) &=& \sigma^{-n+1}(a)\cdot\ldots\cdot a
\end{eqnarray*}
If $n>1$ then the latter have $t$-degree $3n$, hence $K[\alpha(r,-r), r \in \Z] \subset K[t^3]$, so $\alpha(s,t) \notin K[\alpha(r,-r), r \in \Z]$ for all $s,t$ with $s \neq -t$.

\subsubsection{An example of class 3}

Look at the algebra $A$ given by generators $X,Y,Z$ over the field $K$ with relations as follows ($\lambda \neq 0$ in $K$, $c,d \in K$):
\begin{eqnarray*}
XZ &=& \lambda ZX\\
YZ &=& \lambda^{-1}ZY\\
\sqrt{\lambda}YX &=& -(c-Z)(d+Z)\\
\sqrt{\lambda}^{-1}XY &=& -(c-\lambda Z)(d + \lambda Z)
\end{eqnarray*}
Put $t = Z$ and $A_0 = K[t]$, $\sigma{t} = \lambda t$ and 
\begin{eqnarray*}
a = \alpha(-1,1)&=&YX = -\sqrt{\lambda}^{-1}(c- t)(d +t)\\
\sigma(a) = \alpha(1,-1)&=&  -\sqrt{\lambda}^{-1}(c- \lambda t)(d +\lambda t)
\end{eqnarray*}
Put $a'=-\sqrt{\lambda}^{-1}a$ and assume $\lambda \neq 1$, $c \neq d$.  Then we calculate $-\lambda^2 a' + \sigma(a') = (1-\lambda)\lambda(c-d)t +(1-\lambda^2)cd$.  Now
\[K[a,\sigma(a)]=K[a',\sigma(a')]=K[a', -\lambda^2a'+\sigma(a')]=K[a',t]=K[t]\]
Consequently this example is of class 3.

\subsubsection{A new example of general type}\label{19}

Take $A_0=K[t]$, $G=\Z$, $\sigma \in \aut_k A_0$, $\sigma_m = \sigma^m$, $m \in \Z$.  Choose a polynomial $p \in K[t]$ that is irreducible.  Calculating the cocycle-like conditions on $\alpha(n,m)$ yields:
\begin{eqnarray}
\alpha(1,1)\alpha(2,-1) &=& \sigma(\alpha(1,-1))\alpha(1,0)=\sigma(\alpha(1,-1)) \label{13}\\
\alpha(-1,1)\alpha(0,1) &=& \sigma_{-1}(\alpha(1,1))\alpha(-1,2)=\sigma(\alpha(-1,1)) \label{14}\\
\alpha(1,-1)\alpha(0,-1) &=& \sigma(\alpha(-1,-1))\alpha(1,-2)=\alpha(1,-1) \label{15}\\
\alpha(-1,-1)\alpha(-2,1) &=& \sigma_{-1}(\alpha(-1,1))\alpha(-1,0)=\sigma_{-1}(\alpha(-1,1)) \label{16}\\
\alpha(1,-2)\alpha(-1,1) &=& \sigma(\alpha(-2,1))\alpha(1,-1) \label{17}\\
\alpha(1,-1)\alpha(0,1) &=& \sigma(\alpha(-1,1))\alpha(1,0)=\sigma(\alpha(-1,1)) \label{18}
\end{eqnarray}
Assume now that we take $\alpha(-1,1)=p \in K[t]$, then from (\ref{18}) $\sigma(\alpha(1,-1))=\sigma^2(p)$.  From (\ref{13}) and (\ref{14}) it follows then that $\alpha(1,1)=1$.  From (\ref{15}) and (\ref{16}) it follows that $\alpha(-1,-1)=1$.  In that case we obtain again a GWA of degree $1$.\\
However, if we put $\alpha(-1,1) = p \sigma^{-1}(p)$ and $\alpha(1,1)=\sigma(p)$ then we derive
\begin{eqnarray*}
\textup{from (\ref{13}), that: }\ \ \alpha(2,-1)&=&\sigma^2(p)\\
\textup{from (\ref{14}), that: }\ \ \alpha(-1,2)&=&\sigma^{-1}(p)\\
\textup{from (\ref{17}), that: }\ \ \alpha(1,-2)&=&\sigma(p)\\
\                               \ \ \alpha(-2,1)&=&\sigma^{-2}(p)\\  
\textup{from (\ref{15}), that: }\alpha(-1,-1)&=&\sigma^{-1}(p)
\end{eqnarray*}
This leads to the definition of a crystalline graded algebra by putting
\begin{eqnarray*}
\alpha(n,m)&=&\sigma^n(p) \textup{ for } n,m \neq 0 \in \Z \textup{ and } m \neq -n\\
\alpha(n,-n)&=&p\sigma^n(p) \textup{ for } n \neq 0 \in \Z\\
\alpha(0,m)&=&\alpha(m,0)=1 \textup{ for } m \in \Z
\end{eqnarray*}
Observe that $\alpha(1,1) \notin K[\sigma^k(\alpha(n,-n)),k,n \in \Z]$, $\sigma^{-1}(p)=\sigma^{-1}(\alpha(1,1)) \notin K[\alpha(n,m),n,m \in \Z]$.  So now we obtain a crystalline graded ring not in any of the classes we distinguished at the beginning of this section.

\subsubsection{Roll-up examples}

Just like the rolled-up Rees ring appears rather naturally in the theory of filtered rings, in particular in its application to noncommutative geometry (see \cite{LVO} for variations on the Rees ring theme), one may roll-up $\Z$ or ($\Z^n$-) examples to examples for finite cyclic groups.\\
Take $A_0 = \Z_p[X]$ (other grouprings equally possible).  Let $G = \left\langle g\right\rangle$ be a cyclic group of order $n$ and $\sigma \in \aut A_0$, $\sigma^n = I$, write $\sigma_m = \sigma^m$.  Consider an irreducible $p \in k[X]$, $k= \Z_p$, and define
\begin{eqnarray}
\alpha(g^i,g^j) &=& \sigma^i(p) \textup{ if } i,j \neq 0 \modu n \textup{ and } i+j \neq 0 \modu n\\
\alpha(g^i,g^j) &=& p\sigma^i(p) \textup{ if } i+j = 0 \modu n \textup{ but } i,j \neq 0 \modu n\\
\alpha(e,g^j)&=&\alpha(g^j,e)=1 \textup{ for } j=0,1,\ldots, n-1
\end{eqnarray}
It is clear how to define the general construction of a rolled-up $\Z/n\Z$-gradation from a $\Z$-gradation; the example is just a special case of the situation where a crystalline graded ring of the type defined in \ref{19} is rolled up into an example with respect to a cyclic group.\\
Perhaps one of the most interesting generalizations allowed by the definition of crystalline graded rings when compared to GWA is that the automorphisms of $A_0$ do not have to commute, making the techniques available for actions of nonabelian groups, etc.  This may be interesting in very specific cases generalizing the second Weyl algebra for example, e.g. take the twisted tensor product of two first Weyl algebras (see \cite{JLPV} for applications of twisted tensor products to noncommutative manifolds e.g. spherical varieties) and look at the examples of crystalline graded rings constructed this way.  Classical problems, e.g. classification of all simple modules over these, provide interesting problems here (work in progress) where the graded methods (see for example \cite{BVO1} have to be modified but remain valid.)


\begin{thebibliography}{A}

\bibitem[B1]{B1} V. Bavula, Generalized Weyl Algebras and their Representations, Algebra i Analiz 4 (1992), no 1, pp. 75-97.  English translation in St. Petersburg Mat. J. 4 (1993), no 1, pp. 71-92
\bibitem[B2]{B2} V. Bavula, Global Dimension of Generalized Weyl Algebras, CMS Conference Proceedings vol 18 (1996), pp. 81-107.
\bibitem[BVO1]{BVO1} V. Bavula, F. Van Oystaeyen, Simple Holonomic Modules over the Second Weyl Algebra, Adv. Math. 150 (2000), no 1, pp. 80-116.
\bibitem[BVO2]{BVO2} V. Bavula, F. Van Oystaeyen, Krull Dimension of Generalized Weyl Algebras and Iterated Skew Polynomial Rings, J. of Algebra 208 (1998), no 1, pp. 1-34.
\bibitem[BVO3]{BVO3} V. Bavula, F. Van Oystaeyen, The Simple Modules of Certain Generalized Crossed Products, J. of Algebra 194, no 2, pp. 521-566.
\bibitem[J]{J} D.A. Jordan, Krull and Global Dimension of Certain Iterated Skew Polynomial Rings, Contemp. Math 130 (1992), pp. 201-213.
\bibitem[LVO]{LVO} Li Huishi, F. Van Oystaeyen, Zariskian Filtrations, Springer, 1996.
\bibitem[JLPV]{JLPV} P. Jara, J. Lopez, F. Panaite, F. Van Oystaeyen, On Iterated Twisted Tensor Products of Algebras,  arXiv:math.QA/0511280.
\bibitem[N]{N} E. Nauwelaerts, Generalized $2$-Cocycles of Finite Groups and Maximal Orders, J. of Algebra 210 (1998), pp. 225-241.
\bibitem[NaVO1]{NaVO1} C. N\u{a}st\u{a}sescu, F. Van Oystaeyen, Graded Ring Theory, Math. Library vol. 28, North-Holland.
\bibitem[NaVO2]{NaVO2} C. N\u{a}st\u{a}sescu, F. Van Oystaeyen, Methods of Graded Rings, Lect. Notes in Mathematics, Springer Verlag, Berlin 2003.
\bibitem[NeVO]{NeVO} P. Nelis, F. Van Oystaeyen, The Projective Schur Subgroup of the Brauer Group and Root Groups of Finite Groups, J. of Algebra
\bibitem[NVO1]{NVO1} E. Nauwelaerts, F. Van Oystaeyen, Module Characters and Projective Representations of Finite Groups, Proc. London Math. Soc.(3), 62 (1991), pp. 151-166.
137 (1991), no 2, pp. 501-518.
\bibitem[NVO2]{NVO2} E. Nauwelaerts, F. Van Oystaeyen, Finite Generalized Crossed Products over Tame and Maximal Orders, J. of Algebra 101 (1986), pp. 61-68
\bibitem[NVO3]{NVO3} E. Nauwelaerts, F. Van Oystaeyen, Generalized Twisted Group Rings, J. of Algebra 294 (2005), pp. 307-320.
\bibitem[RS]{RS} J.C. Robson, L.W. Small, Liberal Extensions, Proc. London Math. Soc. (3), pp. 87-103.
\bibitem[VO1]{VO1} F. Van Oystaeyen, Azumaya Strongly Graded Rings and Ray Classes, J. of Algebra 103 (1986), no 1, pp. 228-240.
1, pp. 61-68.
\bibitem[VO2]{VO2} F. Van Oystaeyen, Algebraic Geometry for Associative Algebras, M. Dekker Monographs, New York, 2001.

\end{thebibliography}
\end{document}